\documentclass[12pt]{amsart}
\usepackage{latexsym, amssymb, amsmath}

\theoremstyle{plain}

\newtheorem{proposition}{Proposition}
\newtheorem{lemma}{Lemma}
\newtheorem{corollary}{Corollary}

\theoremstyle{definition}
\newtheorem{definition}{Definition}
\newtheorem*{question}{Question}
\newtheorem*{example}{Example}
\newtheorem*{remark}{Remark}

\theoremstyle{remark}

\newcommand{\knot}{\kappa}
\newcommand{\realproj}{{\mathbb R}{\mathbb P}}
\newcommand{\sphere}{\mbox{S}}
\newcommand{\real}{{\mathbb R}}
\newcommand{\complex}{\mathbb C}
\newcommand{\origin}{\boldsymbol 0}

\newcommand{\infplane}{\mbox{H}_\infty}
\newcommand{\xsection}{\S}
\newcommand{\cpoint}{\infty}

\begin{document}

\title{Polynomial knots}

\author{Alan Durfee and Donal O'Shea}

\begin{abstract}
A polynomial knot is a smooth embedding $\knot : \real \to \real^n$
whose components are polynomials.  The case $n = 3$ is of particular
interest.  It is both an object of real algebraic geometry as well as
being an open ended topological knot.  This paper contains basic
results for these knots as well as many examples.
\end{abstract}

\maketitle


\section{Introduction}
\label{introduction}

\begin{definition}
A {\it polynomial knot} is a smooth embedding $\knot : \real \to
\real^n$ whose components are polynomials.  
\end{definition}

A polynomial knot thus has the properties that $\knot'(t) \neq 0$ for
all $t$, and $\knot(s) \neq \knot(t)$ for all $s \neq t$.
Polynomial knots were first investigated by Shastri \cite{Shastri}
in connection with a conjecture of Abhyankar.  He found a
simple representation of the trefoil (overhand) knot:
\begin{align}
\label{Shastri-trefoil}
x(t) &= t^3-3t \\
y(t) &= t^4 -4t^2  \\
z(t) &= t^5-10t. 
\end{align}
This elegant example has remained at the base of the subject.  He
also found equations for the figure-eight knot. We give this knot and
many other examples of equations for topological knots in \xsection
\ref{section-examples}. It is not easy to find these equations.


A topological knot is usually understood to be an embedded circle in
the three sphere.  Polynomial knots are not of this type, but rather a
``long'' or ``open ended'' knot. The one point compactification of a
polynomial knot is a knot in the usual sense.
If the compactification is obtained by using stereographic projection, then
easy arguments show that the compactified polynomial knot is tame, and
also smoothly embedded except possibly at the pole of the sphere where
it has a nicely behaved algebraic singular point.  We also show that a
family of polynomial knots is transverse to a suitably large sphere
about the origin.  This material is in \xsection
\ref{section-basic-results}.

Section \ref{section-approx} is concerned with approximation theorems.
We first show that if $\alpha$ is an embedding of a compact interval
in $\real^n$ and $\beta$ is another map close to both the locus and
derivative of $\alpha$, then $\beta$ is an embedding isotopic to
$\alpha$.  A corollary is that a smooth embedding $\alpha$ on a
compact interval can be approximated by a polynomial $\beta$ which has
the same ``knot type'', and that the same is true if $n = 2$ and the
map $\alpha$ has transverse intersections.  This result is obtained
using the Weierstrass approximation theorem and integration.

Shastri proved that for every topological knot there is a polynomial
knot equivalent to it. The approximation results above fill in some
details of his proof.  These methods also show that a topological knot
can be approximated by a finite Fourier series as well as a rational
function.  A constructive proof of Shastri's result has been provided
by Wright \cite{Wright}.

Akbulut and King \cite{Akbulut-King} prove a related result, namely if
$K \subset \sphere^3 \subset \real^4$ is a compact smooth knot then
there is a real algebraic variety $Z \subset \real^4$ such that $Z
\cap \sphere^3$ is knot-equivalent to $K$.  Also Mond and van Straten
\cite{Mond-vanStraten} show that any knot type can be realized as the
real Milnor fiber of an appropriate space curve singularity.  Note
that Shastri's result guarantees a polynomial parameterization for any
knot type.

In addition to their usual equivalence as topological knots,
polynomial knots have two stronger types of equivalence coming from
the algebraic structure given by the coefficients of their defining
polynomials.  These coefficients form an open subset 
of a high-dimension Euclidean space. We call two polynomial
knots {\em path equivalent} if they are in the same connected
component of 
of this subset.  We show that two path-equivalent polynomial knots are
topologically equivalent.  This argument is not trivial since it is
conceivable that an area of knottedness might slide off to infinity as
the knot moves in a family.

{\em Left-right equivalence} is defined by applying
orientation-preserving linear automorphisms to the domain $\real^1$
and range $\real^n$. A knot obtained by applying such transformations
is easily shown to be both topologically and path equivalent to the
original knot. More generally, one can use polynomial automorphisms.
These results are in \xsection \ref{section-equivalence}.


In \xsection \ref{section-reduced-form} we discuss how the equations
defining a polynomial knot can be reduced to ones of simpler form.  In
fact, a polynomial knot can be reduced, using the standard operations
of linear algebra (adding one row to another and so forth), to one
whose coefficient matrix is in row echelon form.  A stronger operation
is available, that of adding polynomial multiples of one row to
another. This may reduce the degree of the knot. (The {\em degree} of
a polynomial knot is defined as the maximum degree of its component
polynomials.)  Note that the first operation is available in
projective space, whereas the second is not.  These results are used
to find simple equations for a polynomial knot representing a given
topological knot.

In \xsection \ref{section-examples-restrictions} we give many examples,
most of which are due to students with whom we have worked under the
auspices of the National Science Foundation Research Experiences for
Undergraduates (REU) program.
We also show that the
degree of a polynomial representative of a topological knot is bounded
below by expressions in the crossing number, bridge number, and
superbridge number.  The first two of these were found by Lee Rudolph,
and the third by our students. 
These restrictions allows us to list and analyze the topological types
represented by polynomial knots degree by degree. In this task the
first result mentioned is the most useful.

In this section we also discuss curvature and minimal lexicographical
orderings, and we compare the results for polynomial knots with those
for polygonal knots.

In \xsection \ref{section-spaces} we discuss the topological structure
of spaces of polynomial knots.  These results, mostly due to Vassiliev
\cite{Vassiliev}, are augmented by results due to our REU
students. Here again we compare with the situation for polygonal
knots.

Section \ref{section-complex} contains remarks on the situation over
the complex numbers. We discuss the question of whether an embedding
$\complex \to \complex^n$ is rectifiable (left equivalent to a linear
embedding).  We also mention a topic from algebraic geometry, the
standard correspondence of maps $\complex^1 \to \complex^n$ with maps
$\complex[x_1, x_2, \dots, x_n] \to \complex[t]$, as well as the
correspondence of maps which are embeddings to surjective ring maps.

Many other special types of knots have been studied: rational knots
(whose components are rational functions), trigonometric knots (whose
components are finite Fourier series), holonomic knots (knots of the
form $(f, f', f'')$ for some function $f$), polygonal or stick knots
(whose components are piecewise-linear embeddings $\real \to
\real^3$), knots on a cubic lattice, thick knots (made of rope of a
certain radius), and so forth. Polynomial knots are an addition
to these boutique knots; they are attractive because of their
connections with algebraic geometry. 
Knots whose components are rational functions will be the subject of
a later paper.


Many of our results are stated for embeddings $\real \to \real^n$
where $n \geq 2$, which we will also loosely refer to as polynomial
knots.  There is no extra cost in doing this, and some of the
results are interesting in dimensions other than three.

We thank the National Science Foundation for their support of our
Research Experiences for Undergraduates (REU) program 
DMS-9732228 and
DMS-0353700, and we thank the students for their constant
stimulation. We also thank Elmer Rees for helpful conversations, and
the Universities of Miami and Edinburgh for their congenial
atmosphere. Additional support was provided by the Hutchcroft Fund at
Mount Holyoke College.



\section{Basic results}
\label{section-basic-results}

\subsection{Topological knots}
\label{subsection-topological-knots}

We first review some basic concepts of knot theory.  A topological
knot is a continuous embedding $\alpha$ of the circle $\sphere^1$ in
the $n$-sphere $\sphere^n$, where $n \geq 2$. (Usually one takes $n =
3$, but it will be useful to consider the other cases as well.)  
Let $\sphere^1$ and $\sphere^n$ have orientations coming from
the inclusion $\real^k \subset \sphere^k$, where $\real^k$ has its
standard orientation.

Topological knots $\alpha_0$ and $\alpha_1$ are by definition
topologically equivalent (``have the same knot type'') if there are
orientation-preserving homeomorphisms $f$ of $\sphere^1$ and g of
$\sphere^n$ such that
$$\alpha_0  = g \circ \alpha_1 \circ f.$$
Topological equivalence thus preserves the orientations on both the
domain and the range. 

Using stereographic projection we may identify the one point
compactification
of real $n$-space 
with the $n$-sphere.
This map takes the point of compactification 
to the north pole of the $n$-sphere. 
A polynomial embedding $\knot : \real \to \real^n$ extends to a
unique map ${\bar \knot}: \sphere^1 \to \sphere^n$ taking
poles to poles. The map ${\bar \knot}$ is a continuous
embedding, but
not necessarily smooth because there is an algebraic singularity
at the north pole which may or may not be ${\mathcal C}^1$. The details
may be found in \xsection \ref{subsection-stereographic-projection}.  
Even if it has a singularity, a polynomial knot is a knot in the
classical sense:

\begin{proposition}
\label{prop-polynomial-knot-is-tame}
The completed polynomial knot ${\bar \knot}: \sphere^1 \to
\sphere^n$ is tame.
\end{proposition}

Recall that a topological knot is tame if it is topologically
equivalent to a polygonal knot.
Tame knots in dimension three are the subject of classical knot
theory; in dimensions greater than three all tame knots are trivial.
The above proposition follows from the next two lemmas.

\begin{lemma}
The knot ${\bar \knot}$ is piecewise differentiable.
\end{lemma}

By ``piecewise differentiable'' we mean that the domain can be divided
into a finite number of segments such that the map is differentiable
on the interior of each segment and has a limiting tangent line at
each end.  This follows from the fact that the curve has an isolated
algebraic singularity in a neighborhood of the point at infinity.

\begin{lemma}
A piecewise differentiable knot is tame.
\end{lemma}

In fact, Appendix I of \cite{Crowell-Fox} shows that that a
differentiable knot is tame, and the proof easily can be modified: If
$p$ is a cusp point, then pull apart the branches of the cusp so that
they are no longer tangent to each other, then make each cusp point
point the vertex of a cone.


We can also compactify
$\real^n$ by projective $n$-space $\realproj^n$. The space
$\realproj^n - \real^n$ is $ \infplane^n$, the hyperplane at infinity.
(Note that $\infplane^1$ is a point.)
A polynomial map $\knot : \real \to \real^n$
extends to a unique map ${\tilde \knot}: \realproj^1 \to
\realproj^n$ taking $\infplane^1$ to a point on $\infplane^n$. The
projection $\realproj^n \to \sphere^n$ collapsing the hyperplane
$\infplane^n$ to the north pole of $\sphere^n$
takes $\tilde{\knot}$ to ${\bar \knot}$. 
Projective
compactifications will be discussed more fully in the forthcoming paper
\cite{Durfee-O'Shea-rational}.

\subsection{Stereographic projection and the singularity at infinity}
\label{subsection-stereographic-projection}

In this section we use stereographic projection to give an analytic
chart for $\real^n$ at infinity, then use this to describe the
singularity of a polynomial knot at this point.

Now let $\sphere^n \subset \real^{n+1}$ be the unit sphere centered at
the origin, $n \geq 1$.
The map  
$$\sigma': \real^n \times \{0\}  
\longrightarrow \sphere^n $$
with formula
$$({\bf x}, 0) \longrightarrow \left(
\frac{2{\bf x}}  {{| {\bf x} |}^2 + 1}, \ \frac{ {| {\bf x} |}^2
- 1}  { {| {\bf x} |}^2 +1} \right). $$ 
yields the stereographic projection map 
$$ \sigma: \real^n \cup \{\infty\} \longrightarrow \sphere^n, $$
taking the point of compactification $\infty$ to the
north pole of $\sphere^n$.  
We use this map to provide an analytic chart for $\real^n$ in a
neighborhood of $\infty$.

Let $\knot$ be a polynomial map of degree $d > 0$. 
The compactification $\bar{\knot}$ of  $\knot$ under stereographic
projection  is the 
closure of the composition $\sigma\circ\knot$.
As $t$ tends to $\pm \infty$, the composition $\sigma\circ\knot(t)$
tends to the north pole.
Suppose that 
\begin{equation} \label{knot-equation}
\knot(t) = {\bf a}_dt^d + {\bf a}_{d-1}t^{d-1} + \ldots + {\bf
a}_1t +{\bf a}_0
\end{equation}
where ${\bf a}_d, \ {\bf a}_{d-1}, \ldots, {\bf a}_0 \in
\real^n$ and ${\bf a}_d \ne \origin$. 
The standard coordinates on $\real^n$ near the origin give local
coordinates on $\sphere^n$ near the north pole. In terms of these, the
map $\bar \knot$ is given by
\begin{align}
t \mapsto \frac{2\knot(t)}{| \knot(t)|^2+1}.
\end{align}
Setting $s = 1/t$, the above expression becomes
$$
\frac
{2 ( {\bf a}_ds^d + {\bf a}_{d-1}s^{d+1} + \ldots + {\bf a}_0s^{2d} ) }
{ {\bf a}_d\cdot {\bf a}_d + 2{\bf
a}_d\cdot{\bf a}_{d-1} s + \ldots + ({\bf a}_0\cdot{\bf a}_0
  +1)s^{2d}}
$$
which simplifies to
$$
{2{\bf a}_d\over |{\bf a}_d|^2} s^d + [h].
$$ 
(Here, and below, $[h]$ denotes higher order terms, and $[l]$ lower
order terms.)  The computations above establish the following result.

\begin{proposition}
\label{prop-singularity-at-infinity}
Let $n \geq 1$, let $\knot: \real \rightarrow \real^n$ be a
polynomial map of degree $d$ and let $\sigma: \real^n \cup \{\infty\}
\rightarrow \sphere^n$ be stereographic projection.  The
compactification ${\bar \knot} = \sigma\circ\knot$ is a curve in
$\sphere^n$ with an analytic singularity of order $d$ in local
coordinates at the north pole.
\end{proposition}

For example, the Shastri trefoil $\knot(t) = (t^3-3t, t^4-4t^2,
t^5-10t)$ has singularity at the pole with expansion 
$2(s^7, s^6, s^5) + [h]$.
More generally, a polynomial knot of the form
$\knot(t) = (t^p, t^q, t^r) + [l]$
with $p \le q \le r$
has local analytic expansion 
$2(s^{2r-p}, s^{2r-q}, s^r) + [h]$ 
at the north pole. 
In particular a polynomial knot with $d > 1$
will always have an algebraic singularity at infinity. It may however
be  ${\mathcal C}^1$ there.

\subsection{Transversality}
\label{subsection-transverse}

In this section we show that polynomial maps are well-behaved at
infinity, in fact transverse to all sufficiently large spheres
about the origin.

\begin{proposition}
\label{proposition-transverse} 
(a). If $n \geq 1$ and $\alpha: \real \rightarrow \real^n$ is a
polynomial map, then 
as $t \to + \infty$
the angle between the vectors $\alpha(t)$ and
$\alpha'(t)$ approaches 0, 
and as $t \to - \infty$ this angle approaches $\pi$.
In particular, there exists $r_0 \gg 0$ such that the $(n-1)$-sphere
of radius $|\alpha(t)|$ about the origin is transverse to
$\alpha'(t)$ for all $|t| \geq r_0$.

\noindent (b). Suppose that $\{ \alpha_u \}$ is a family of maps
depending continuously on a parameter $0 \le u \le 1$. If the degree
of $\alpha_u$ is constant (independent of $u$), then there is
an $r_0$ such that Part (a) is true for all $\alpha_u$. 
\end{proposition}



\begin{proof}  

\noindent (a). Suppose that $\alpha(t)$ is as in Equation
\ref{knot-equation}.  
The cosine of the angle between the vectors $\alpha(t)$ and $\alpha'(t)$
is 
$$
\frac{\alpha(t)}{|\alpha(t)|} \cdot \frac{\alpha'(t)}{|\alpha'(t)|} 
$$
which becomes
$$
\pm \frac
{ {\bf a}_d + \frac{1}{t} {\bf a}_{d-1} + \cdots + \frac{1}{t^d}{\bf a}_0}
{| {\bf a}_d + \frac{1}{t} {\bf a}_{d-1} + \cdots + \frac{1}{t^d}{\bf a}_0|}
\cdot
\frac
{d {\bf a}_d + \frac{d-1}{t} {\bf a}_{d-1} + \cdots 
  + \frac{1}{t^{d-1}}{\bf a}_0}
{| {\bf a}_d + \frac{d-1}{t} {\bf a}_{d-1} + \cdots 
  + \frac{1}{t^{d-1}}{\bf a}_0 |} \ .
$$
As $t \to \pm \infty$ this expression approaches
$$
\pm \frac
{ {\bf a}_d \cdot d {\bf a}_d }
{| {\bf a}_d \cdot d {\bf a}_d |}
= \pm 1 .
$$

\noindent (b).
The coefficients of $\alpha_u(t)$ are functions of  $u$.  
Since ${\bf a}_d \neq 0$ 
the denominators can be bounded below by a nonzero constant.  The
above argument then continues to hold.
\end{proof}



\section{Polynomial approximations of smooth knots}

\label{section-approx}

\subsection{Approximation}   

Let $M \subset \real$ be a compact connected interval.  In this
section we show that given a smooth embedding $\alpha: M \to \real^n$,
and another map $\beta: M \to \real^n$ close to it (both the function
and its derivative), then $\beta$ is an embedding which is isotopic to
$\alpha$.

\begin{lemma}
\label{lemma-inequalities}
Suppose that $\gamma: M \to \real^n$ is ${{\mathcal C}^2}$, where $n \geq 1$. 
\newline (a).  If $\gamma'(t) \neq 0$ for all $t \in M$, then there
are positive numbers $r$ and $\delta$ such that $|s-t| \leq r$ implies
that  
$ |\gamma(s) - \gamma(t) | \geq \delta |s-t|$.
\newline (b). If $\gamma(s) \neq \gamma(t)$ for all $s
\neq t \in M$, then for each $r>0$ there is a $\delta$ such that
$|s-t| \geq r$ implies that $|\gamma(s) -\gamma(t)| \geq \delta$.
\end{lemma}

\begin{proof}
(a). Since $M$ is compact and $\gamma'(t) \neq 0$ for all $t \in M$, there
is a $\delta > 0$ such that $|\gamma'(t)| \geq 2\delta$ for all $t \in
M$. 
Let $m$ be the maximum value of $(1/2) |\gamma''(t)|$ for $t \in M$.
If $m = 0$ then $\gamma$ is linear and the result is clear.  
Otherwise let
$r = \delta / m$.
By Taylor's theorem, 
for each $s,t \in M$ there exists $\tau$ between $s$
and $t$ such that
\begin{align*}
|\gamma(s) - \gamma(t)| & = | \gamma'(t)(s-t) + (1/2)\gamma''(\tau)
(s-t)^2| \\
& \geq |s-t| \Bigl| |\gamma'(t)| - | (1/2) \gamma''(\tau) (s-t) |
\Bigr| \\
& \geq \delta |s-t|.
\end{align*}
The last line follows since $|s-t| \leq r$ implies that 
$(1/2) |\gamma''(\tau)| |s-t| \leq mr = \delta$.
\newline (b). This follows since the set 
$K = \{ (s,t) \in M \times M : |s-t| \geq r \} $ is compact
so the function $|\gamma(s) - \gamma(t)|$ on $K$ is bounded below.
\end{proof}

\begin{definition}
Let $M \subset \real$ 
and let
$\alpha, \ \beta: M \to \real^n$ be ${\mathcal C}^1$ maps. The maps
$\alpha$ and $\beta$ are $\epsilon$-{\it close} if $|\alpha(t) -
\beta(t)| < \epsilon$ and $|\alpha'(t) - \beta'(t)| < \epsilon$, for
all $t \in M$.
\end{definition}

\begin{proposition}

\label{approx-prop}
(a). Suppose that $n \geq 1$, and that $M$ is a compact connected interval.
Let $\alpha: M \to \real^n$ be a $\mathcal{C}^2$ embedding and let
$\beta: M \to \real^n$ be a ${\mathcal C}^2$ map $\epsilon$-close to
$\alpha$, where $\epsilon > 0$ is a suitably small
number.
Then $\beta$ is an embedding, and $\beta$ is 
isotopic to $\alpha$.
\newline (b). If $n = 2$, and the word ``embedding'' is replaced by
``has transversal self-intersections'', then the same proposition holds.
\end{proposition}

\noindent The number $\epsilon$ will be determined in the course of
the proof. 

\begin{proof}
Let
$$\gamma_u(t) = u\beta(t) + (1-u) \alpha(t)$$ where $0 \leq u \leq 1$.
Note that $\gamma_0 = \alpha$ and $\gamma_1 = \beta$. Fix $u$. We will
show that $\gamma_u$ is an embedding. The family
$\gamma_u$ thus provides an isotopy of $\alpha$ to $\beta$.

First we show that $\gamma'_u(t) \neq 0$, for all $t \in M$:
\begin{align}
|\gamma'_u(t)| &=
|u\beta'(t) + (1-u) \alpha'(t) |\nonumber \\
&\geq
\Bigl| | \alpha'(t)| - u|\beta'(t) - \alpha'(t)| \Bigr|\nonumber \\
&\geq \Bigl| | \alpha'(t)| - |\beta'(t) - \alpha'(t)|
\Bigr|.   \nonumber
\end{align}
Since $M$ is compact and $\alpha' \neq 0$ for $t \in M$, there is a
$\delta_1 > 0$ such that $|\alpha'(t)| \geq \delta_1$. If $|\beta'(t)
- \alpha'(t)| \leq \delta_1/2$ then the above expression is at least
$\delta_1/2$ and hence nonzero.

Next we show that $s \neq t$ implies that $\gamma_u(s) \neq
\gamma_u(t)$. Since $\gamma'_u(t) \neq 0$, Lemma \ref{lemma-inequalities}(a)
implies that there is an
$r>0$ such that such that this is true if $|s-t| \leq r$.
Now suppose that $|s-t| \geq r$. Then
\begin{align}
\Bigl| \gamma_u(s) - \gamma_u(t) \Bigr| &=
\Bigl| [u\beta(s) + (1-u)\alpha(s)] -[ u\beta(t) + (1-u) \alpha(t)] \Bigr|
\nonumber \\
&= \biggl| \Bigl[\alpha(s) - \alpha(t) \Bigr]
+ u\Bigl[(\beta(s) -\alpha(s) ) + (\alpha(t) - \beta(t)) \Bigr] \biggr|\nonumber \\
&\geq \biggl| | \alpha(s) - \alpha(t)| - u\Bigl[ |\beta(s) -
\alpha(s)| + |\alpha(t)-\beta(t)| \Bigr] \biggr|.\nonumber
\end{align}
By Lemma \ref{lemma-inequalities}(b), there is a $\delta_2 > 0$ such
that $|s-t| \geq r$ implies $|\alpha(s) - \alpha(t)| \geq
\delta_2$. If $ |\beta(t) - \alpha(t)| \leq \delta_2/4$, then the
above expression is nonzero.
If $0 < \epsilon \leq \mbox{min} \{\delta_1/2, \delta_2/4\}$
the then statement (a) follows.  The proof of statement (b) is
essentially the same.
\end{proof}

\subsection{Polynomial knots} 
\label{poly-knot-section}

We now apply the above results to fill in some details in Shastri's
result that for every topological knot there is a polynomial knot of
the same knot type. The polynomial knot constructed here will
not be $\epsilon$-close to the original knot; in fact the original is
a compact subset of $\real^n$ whereas the polynomial knot will pass
through the point at infinity.

An alternate proof of this approximation theorem is due to
Wright \cite{Wright}.  The method here is first to replace the $z$-coordinate
of the knot by an interpolating polynomial which goes under and over
at the same times.  The result is then turned on its side so that the
$x$-coordinate is up, and this is also approximated by an
interpolating polynomial.  This is repeated with the $y$-coordinate.
The final result is a polynomial knot of the same knot type as the
original. This method is algorithmic, and can be turned into a computer
program. 

In addition, Mui \cite{Mui, Mui-thesis} describes a method to find 
height functions of low degree for a given knot projection.

\begin{proposition} 
\label{poly-approx-thm}
(a). If $n \geq 1$ and $\alpha : [a,b] \to \real^n$ is a ${\mathcal C}^2$ 
embedding, then there is a polynomial embedding 
$\beta: [a, b] \to \real^n$ 
isotopic to $\alpha$. Furthermore the maps $\alpha$ and $\beta$ can be
taken $\epsilon$-close for any small $\epsilon > 0$.
\newline (b). If $n = 2$, and the word ``embedding'' is replaced by
``has transversal self-intersections'', then the same proposition holds.
\end{proposition}


\begin{proof}
Let $\xi$ be a small positive number.  By the
Weierstrass approximation theorem there is a polynomial $\gamma$ such
that $|\gamma(t) - \alpha'(t)| < \xi$ for all $t \in [a,b]$.
If
$$\beta(t) = \int_a^t \gamma(s) ds + \alpha(a)$$
then $\beta'(t) = \gamma(t)$, and
\begin{align*}
|\beta(t) - \alpha(t)| &= \Bigl| \Bigl[ \int_a^t \gamma(s) ds +
\alpha(a) \Bigr] - \Bigl[ \int_a^t \alpha'(s) ds + \alpha(a)
\Bigr] \Bigr| \\
& = \Bigl| \int_a^t \gamma(s) ds - \int_a^t \alpha'(s) ds \Bigr| \\
& \leq (b-a) \sup_{s \in [a,b]} | \gamma(s) - \alpha' (s) | \\
& \leq (b - a) \xi.  
\end{align*}
If $\xi$ is small enough so that $\epsilon = \mbox{max} \{ \xi,
(b-a)\xi \}$ satisfies the hypotheses of Proposition
\ref{approx-prop}, then part (a) follows.  The proof of part (b) is
similar. 
\end{proof}


The above proposition fills in a detail in the following result due to
Shastri:

\begin{proposition} For every ${\mathcal C}^2$ knot there is a
  polynomial knot of the same knot type.
\end{proposition}

In fact, let 
$$
\alpha(t) = (\alpha_1(t), \alpha_2(t), \alpha_3(t)): 
\realproj^1 \to \real^3
$$ 
be a ${\mathcal C}^2$ embedding with the required knot type.
(We use the real projective line  $\realproj^1$ rather than the circle
$\sphere^1$ in order to have the inclusion $\real \subset \realproj^1$.)
Without loss of generality we may assume that the map
$$
t \to (\alpha_1(t), \alpha_2(t)): \realproj^1 \to \real^2
$$ has only double points.  Choose an interval $[a,b] \subset \real^1
\subset \realproj^1$ containing these double points.  Find a
polynomial function $\beta_3: \real^1 \to \real^3$ which represents
the height function $\alpha_3$.  There is $c < a< b < d$ with $c,d \in
\real$ such that $\beta_3$ is monotonic on $(-\infty, c)$ and $(d,
+\infty)$.  
By Proposition \ref{poly-approx-thm}(b) there is a
polynomial map $(\beta_1, \beta_2): [c, d] \to \real^2$ which has
transverse intersections and which is isotopic to the map $(\alpha_1,
\alpha_2)$.  The spherical completion of 
$$
\beta = (\beta_1, \beta_2, \beta_3): \real^1 \to \real^3
$$ 
has
the same knot type as $\alpha$.

If the polynomial approximation $(\beta_1, \beta_2)$ to the plane
projection had been chosen before the height function $\beta_3$, then
the result may not have the correct knot type since $(\beta_1,
\beta_2)$ may have additional intersections outside of the interval of
approximation.

\subsection{Trigonometric knots}

Suppose that $n \geq 1$. 
A finite Fourier series $\beta: [0, 2\pi] \to \real^n$ is a map of
the form
$$ \beta(t) = \sum_{k = 0}^m {\bf a}_k \sin kt + {\bf b}_k \cos
kt$$
where ${\bf a}_k$ and ${\bf b}_k$ are real vectors of dimension $n$. 

\begin{proposition}
\label{prop-trig-approx}
Let $\epsilon > 0$.
If $\alpha : \real \to \real^n$ is a ${\mathcal C}^2$
embedding, and if $\alpha$ is periodic with period $2\pi$, then there
is a finite Fourier series which is an $\epsilon$-close embedding
isotopic to $\alpha$.
\end{proposition}

The proof is similar to the corresponding proof for polynomials, using
the theorem that $\alpha$ has a Fourier series whose partial sums
converge uniformly to $\alpha$.

\subsection{Rational knots}

\begin{proposition}
Let $\epsilon > 0$.
If $\alpha : \realproj^1 \to \real^n$ is a ${\mathcal C}^2$ embedding,
then there is a rational function which is an $\epsilon$-close embedding 
isotopic to $\alpha$.
\end{proposition}

By Section \ref{poly-knot-section} there is a polynomial knot of the
same knot type as $\alpha$, but this approximation may not be
$\epsilon$-close. 

\begin{proof}
By Proposition \ref{prop-trig-approx} there is a finite Fourier series
$\beta(t)$ which is an embedding  isotopic to
$\alpha$. The standard substitution
$$t = \frac{1}{2} \arctan x$$ changes this into a rational function.
(For details see Clark \cite{Clark}.)
\end{proof}

The rational functions produced by this method appear to have
needlessly high degree, so it would be nice to have a direct proof of
this proposition.



\section{Equivalence of polynomial knots}
\label{section-equivalence}


The polynomial knots in this section are of the form $\knot: \real \to
\real^n$ for $n \geq 2$.  
We give three different definitions of equivalence between two
polynomial knots $\knot_0$ and $\knot_1$ and describe connections among them. 

The three types of equivalence are: (1) have the same topological knot
type, (2) are in the same connected component of the coefficient
space, and (3) if transformations of the range and domain take one to
the other.  If they are equivalent in either the second or third
sense, then they are equivalent in the first sense. Whether the
converse is true is not as clear, as is the relation between the
last two.


\begin{definition}
Two polynomial knots $\knot_0$ and $\knot_1$ are {\em topologically
equivalent} if there are
orientation-preserving homeomorphisms $f$ of $\sphere^1$ and g of
$\sphere^n$ such that
$$\bar{\knot}_0 = g \circ \bar{\knot}_1 \circ f.$$
\end{definition}

They are thus topologically equivalent if their
spherical completions are equivalent as topological knots.  Recall
that polynomial knots are tame, so the above is the usual knot
equivalence.  Topologically equivalent polynomial knots may have
different degrees; for example one can replace the parameter $t$ by
$t^k + t$ where $k$ is an odd integer.

The next definition of equivalence makes use of the parameter space of
these knots.  Let ${\mathcal M}^n_d$ be the space of all polynomial
maps $\real \to \real^n$ of degree $d$.  This space is isomorphic to a
proper subset of $(\real^{d+1})^n$ and inherits its topology.  (The
inclusion is proper since ${\bf a}_d \neq \origin$.)  Let $\Sigma^n_d
\subset {\mathcal M}^n_d$ be the space of maps $\alpha$ with either
double points ($r \neq s$ such that $\alpha(r) = \alpha(s)$) or
critical points ($t$ such that $\alpha'(t) = 0$).  This is a closed
subset of ${\mathcal M}^n_d$. 
Let ${\mathcal K}^n_d ={\mathcal M}^n_d
- \Sigma^n_d $; this is the parameter space of polynomial
knots of degree $d$.


\begin{definition}
\label{definition-path-equivalence}
Two polynomial knots $\knot_0$ and $\knot_1$ of degree $d$ are {\em
path equivalent} if they are in the same connected component of
${\mathcal K}^n_d$.
\end{definition}

\begin{proposition}
\label{path-implies-topological} If two polynomial knots are path
equivalent, then they are topologically equivalent.
\end{proposition}

\begin{proof}
Let $\knot_0$ and $\knot_1$ be knots in the same path component of
${\mathcal K}^n_d$, and let $\{\knot_u$, $0 \leq u \leq 1\}$ be a
family of knots connecting them.  Since ${\mathcal K}^n_d$ is open we
may assume without loss of generality that this family is smooth in
$u$.  Let $\Theta: \real \times I \to \real^n \times I$ be defined by
$\Theta(t, u) = (\knot_u(t), u)$. (Here $I$ is the unit interval.)
Let $K$ be the image of $\Theta$.  This is a smooth submanifold of
$\real^n \times I$, since each $\knot_u$ has no singularities.  We
identify the one point compactification $\real^n \cup \{ \cpoint \}$
with $\sphere^n$. Let $\text{N}$ be the north pole of $\sphere^n$. Let
$\bar{K} \subset \sphere^n \times I$ be the closure of $K$.  By
Proposition \ref{proposition-transverse} the stratification $\sphere^3
\times I - \bar{K} \supset K \supset \{\text{N}\} \times I$
satisfies the Whitney conditions.  Furthermore the map $\pi: \sphere^3
\times I \to I$ is proper, and for each stratum $S$ the map $\pi|S : S
\to I$ is a submersion. By Thom's first isotopy lemma (\cite{Thom},
\cite{Mather}), this bundle is locally trivial (${\mathcal C}^0$, not
necessarily ${\mathcal C}^\infty$).  Thus the topological knots
$\bar{\knot}_0$ and $\bar{\knot}_1$ are equivalent (preserving
orientations).
\end{proof}

\begin{question}
Is the converse true?  In other words, can two
topologically-equivalent knots be connected by a path of polynomial
knots?  This is an attractive question whose answer is probably
negative, though no examples are known. For example, in dimensions
greater than three all topological knots are unknotted, but there may
be more that one component of unknotted polynomial knots.
\end{question}

\begin{remark}
Suppose that $\gamma_1: \real \to \real^3$ is a open-ended topological
trefoil with the property that $\gamma(t) = (t,0,0)$ for sufficiently
large $|t|$. Let $\gamma_u = \gamma_1 + (1/u - 1,0,0)$ for $0 < u \leq
1$ be this knot translated to the right, and let $\gamma_0(t) =
(t,0,0)$. In this family the knot $\gamma_1$ is translated to infinity and
becomes unknotted.
Here the stratification above does not satisfy the
Whitney conditions, so the proposition does not hold.
\end{remark}

The next type of equivalence combines left equivalence (algebraic
transformations of the range) and right equivalence (algebraic
transformations of the domain).  The first moves the knot around, and
the second reparameterizes it.

The simplest type of transformation that can be used is an affine
automorphism, a map of the form $
T({\bf x}) + {\bf c}$ for 
${\bf x}
\in \real^n$, where $T: \real^n \to \real^n$ is a
bijective linear transformation and \, ${\bf c}$ \, a constant.
More generally,  $T$ can be a map of the form
$T({\bf x}) = (T_1 
({\bf x}), T_2 ({\bf x}), \dots T_n({\bf x}))$ where each $T_i$ is a
polynomial in ${\bf x}$;
this map is a {\em polynomial automorphism} it if has an inverse of
this form.  
For example, the map $(x,y) \mapsto (x, y + x^2)$ of $\real^2$ is a
polynomial automorphism.
If $T_l$ and $T_r$ are polynomial automorphisms and $\knot$ is a
polynomial knot, then so is $T_l \circ \knot \circ T_r$.
Note that a bijective polynomial map with non-zero differential at
each point of $\real^n$ need not be a polynomial automorphism, for
example the map $t \mapsto t^k + t$ of the real numbers with $k$ odd.
Also note that a polynomial automorphism may not extend to projective
space; this makes the theory of affine knots different from the theory
of projective knots.


\begin{definition}
Two polynomial knots $\knot_0, \knot_1: \real \to \real^n$ are {\em
left-right (LR) equivalent} if there are orientation-preserving
polynomial automorphisms $T_r$ of $\real^1$ and $T_l$ of $\real^n$
such that
$$\knot_0 = T_l \circ \knot_1 \circ T_r.$$
\end{definition}

Note also that an affine automorphism of $\real^n$ preserves the
degree of a polynomial knot, whereas a polynomial automorphism may not
(for example the transformation $(x, y) \mapsto (x, y+x^2)$ applied to
the trivial knot $t \mapsto (t,0)$).

Of course there can be many different variations on this type of
equivalence: equivalence using linear transformations, equivalence
using transformations which are not necessarily orientable, left and
right equivalence separately, and so forth. We specify these as
needed.

\begin{proposition}
If two polynomial knots are LR-equivalent, then they are
topologically equivalent.
\end{proposition}

\begin{proof} The map $\bar{T}_l: \sphere^3 \to \sphere^3$, where the
  bar denotes completion, takes $\bar{\knot}_0$ to $\bar{\knot}_1
  \circ (\bar{T}_r)^{-1}$.  The latter is a reparameterization of
  $\bar{\knot}_1$ and hence of the same knot type.
\end{proof}

\begin{proposition}
If two polynomial knots are LR-equivalent by (orientation-preserving)
affine transformations, then they are path equivalent.
\end{proposition}

This follows since the group of such transformations is connected.

\begin{question}
In general, the relation between path equivalence and LR-equivalence
is not clear.  If two polynomial knots of the same degree are
LR-equivalent by polynomial transformations, are they path equivalent?
This is true if the transformations of the domain and range are
tame. (Tame transformations are by definition a composition of maps
which add a polynomial multiple of one row to another row. Not all
polynomial automorphisms are tame; see \cite{Shestakov-Umirbaev}.)
Conversely, if two polynomial knots are path equivalent, are they
LR-equivalent? This is probably not true, though no examples are
known.  For the complex case of these two questions see \xsection
\ref{section-complex}.
\end{question}


\section{Examples and restrictions}
\label{section-examples-restrictions}

In this section we examine the following questions:
\begin{enumerate}
\item Given a topological knot, find a polynomial knot such that
its spherical compactification has  the same knot type.
\item Given a topological knot, find a lower bound for the degree
of a polynomial representation of this knot.
\item Find a polynomial representative of this degree.
\end{enumerate}
We also examine the above questions using the lexicographic ordering
on the degrees of $(x(t), y(t), z(t))$. We compare all these results
with the corresponding situation for stick knots.  Finally, we relate
the total curvature of the knot with its degree.

\subsection{Examples}
\label{section-examples}
Let us start with the first question.  Although every topological knot
has a polynomial representation (\xsection \ref{section-approx}), less
is known about constructing polynomial representatives in low
degree. In this section we give many examples which have low degree
equations. (See \cite{Toman} for a selected list.) Some have elegant
equations, others messy. The latter arise by constantly
adjusting coefficients until the required knot type is obtained.

\subsubsection{The trefoil knot}

Shastri's equations have degree (3,4,5). These are the simplest set of
equations for a nontrivial knot. 

\begin{align}
x(t) &= t^3-3t \nonumber \\
y(t) &= t^4 -4t^2 \nonumber \\
z(t) &= t^5-10t. \nonumber
\end{align}

\subsubsection{The figure-eight knot}

 Shastri also found the following equations of degree $(3,5,7)$:

\begin{align}
x(t) &= t^3-3t \nonumber \\
y(t) &= t(t^2-1)(t^2-4) \nonumber \\
z(t) &= t^7 -42 t. \nonumber
\end{align}

Equations of the same degree but with just two terms were found
by Brown \cite{Brown}:

\begin{align}
x(t) & = t^3 - 5t  \nonumber \\
y(t) & =  t^5 - 28t  \nonumber \\
z(t) & =  t^7 - 32 t^3. \nonumber
\end{align}

The following equations of degree
$(4,5,6)$  were found by McFeron \cite{McFeron-fig8}:


\begin{align*}
x(t) &= -t^4 + 2.279283653 t^3 + 5t^2 -8.63068748t+0.35140383   \\
y(t) &= t^5-5t^3+4t   \\
z(t) &= (t+2.06) \cdot (t+1.916737670) \cdot (t+0.2122155248)  \\
& \ \ \ \ \cdot (t-1.379221313) \cdot (t-2.05) \cdot (t+10). 
\end{align*}

\subsubsection{Knots with five crossings}

\noindent The $5_1$ knot (the torus knot of type (2,5)):
\begin{itemize}
\item Auerbach \cite{Auerbach} (degree (4,5,7)):
\begin{align}
x(t) & = (t^2 + 5t + 4)(t^2 - 7t + 10)   \nonumber \\
y(t) & =  (t^2-5.2t)(t^2 - 9)(t+4.7) \nonumber \\
z(t) & = t^7 - 3.90763t^6-25.8835t^5+ 83.4739t^4 + 176.691t^3
\nonumber \\
     & \ \ \ - 364.064t^2 - 331.888t + 321.285. \nonumber
\end{align}
\item REU 1998  (degree (5,6,7):
\begin{align}
x & = 1000t^5 - 541000t^3 + 44100000t \nonumber \\
y & = 100t^6 - 80400t^4 + 17126400t^2 - 792985600 \nonumber \\
z & = t^7 -  \frac{13433}{16} \, t^5 + \frac{783769}{4} \, t^3 -
9363600 \, t. \nonumber 
\end{align}
\item Brown \cite{Brown} (degree (4,5,7)):
\begin{align}
x(t) & =  t^4 - 24t^2  \nonumber \\
y(t) & =  t^5 - 36t^3 + 260t  \nonumber \\
z(t) & =  t^7 - 31t^5 + 168t^3 + 560t. \nonumber
\end{align}
\end{itemize}

\noindent The $5_2$ knot (REU 2006), equations of degree (4, 5, 9):
\begin{eqnarray*}x(t)&=&t^4-12t^2\\
y(t)&=&t(t^2-4)(t^2-11)\\
z(t)&=&(t-0.5)(t-10)(t-12.5)(t-20)\cdot\\
&&(t-28)(t-29)(t-40)(t-50.2)(t-50.6).
\end{eqnarray*}

\subsubsection{Knots with six crossings}
The $6_1$ knot 
(REU 2006), equations of degree (4, 7, 11):

\begin{eqnarray*}x(t)&=&(t^2-0.5)(t^2-15.3)\\
y(t)&=&t(t^2-16)(t^2-7)(t^2-4.5)\\
z(t)&=&t(t^2-15.2)(t^2-9)(t^2-6.25)(t^2-1)(t^2-0.25).
\end{eqnarray*}

\noindent The $6_2$ knot \cite{Brown}, equations of degree (4, 5, 11):

\begin{eqnarray*}x(t)&=&t^4-12t^2\\
y(t)&=&t(t^2-4)(t^2-11)\\
z(t)&=&t(t^2-1)(t^2-9)(t^2-\frac{49}{16})(t^2-\frac{169}{16})(t^2-\frac{100}{9}). 
\end{eqnarray*}

\noindent The $6_3$ knot (Curry, REU 1998), equations of degree
(3, 8, 10)):

\begin{align*}
x(t) &= t^3 - 100 t \\
y(t) &= t \cdot (t+4) \cdot (t-6) \cdot (t-8) \cdot (t-9) \cdot (t+9)
\cdot (t-11) \cdot (t+11) \\ 
z(t) &= (t+1) \cdot (t^2-4) \cdot (t-6) \cdot (t+7) \cdot (t+9.2) \cdot (t-9.5) \\
    & \ \ \ \ \ \ \ \ \ \ \ \ \  \cdot (t+10) \cdot (t-10.5) \cdot (t+11.5).
\end{align*}


\subsubsection{Other knots}
\label{section-other-examples}
Many other examples of polynomial knots are known: Auerbach
\cite{Auerbach} found equations of degree (3,10,11) for the (2,7)
torus knot, and of degree (3,16,17) for the (2,9) torus knot.  Curry
(REU 1998) found the following equations of degree (7, 6, 7) for the
(4,3) torus knot (the knot $8_{19}$):
\begin{align*}
x(t) &= t(t^2 - 12^2)(t^2 - 28^2)(t^2 - 30^2) \\
y(t) &= (t^2 - 6^2)(t^2 - 23^2)(t^2 - 29^2) \\
z(t) &= -t(t^2 - 8^2)(t^2 - 12.2^2)(t^2 - 29^2).
\end{align*}
Equations for the other knots with the same projection as $8_{19}$ can
be found in \cite{List}.
Mishra and others in a series of papers \cite{Ranjan-Shukla, Mishra99,
Mishra00, Madeti-Mishra-2006a, Madeti-Mishra-2006b} find polynomial
representations of torus and two-bridge knots.  Wright \cite{Wright}
gives an algorithmic method for finding a polynomial representation of
any knot (see \xsection \ref{poly-knot-section}).


\subsection{Reduced forms}
\label{section-reduced-form}

We now show that polynomial knots can be reduced to ones of simpler
form.  This procedure has been used extensively when finding equations
for a specific knot (\xsection \ref{section-examples-restrictions})
and in fact can be used to reduce these equations even more. It
has also been used in \cite{Mishra00}. 
Our reductions will use polynomial automorphisms of the range;
we refer to this as ``left equivalence''. 
Also recall that path equivalence was defined in \xsection
\ref{section-equivalence}. 
Let
$$
\knot(t) = {\bf a}_dt^d + {\bf a}_{d-1}t^{d-1} + \ldots + {\bf
a}_1t +{\bf a}_0
$$ 
be a polynomial knot of degree $d$,
where ${\bf a}_d, {\bf a}_{d-1}, \ldots, {\bf a}_0 \in
\real^n$ for $n \geq 1$ are vectors, and ${\bf a}_d \ne \origin$.
The next result follows by a translation of the range.

\medskip

\noindent {\bf Reduction 1}. The knot $\knot$ is both left and path
equivalent to one with ${\bf a}_0 = {\bf 0}$.

\medskip

In what follows we assume that ${\bf a}_0 = {\bf 0}$.
We define the {\em coefficient matrix} of $\knot$ to be 
$$
A = [ {\bf a}_d \ {\bf a}_{d-1} \cdots {\bf a}_1 ]
$$
where the ${\bf a}_i$ are regarded as column vectors.
The next two reductions are based on the simple observation that
$$ 
L \circ \knot(t) = LA \, {\bf t}, 
$$
where
${\bf t}$ is the column vector $( t^d, t^{d-1}, \dots, t)$
and $L$ is a linear transformation of $\real^n$.
In other words, linear transformations of the range are the same as
left multiplication of the coefficient matrix. 

\medskip

\noindent {\bf Reduction 2}. The knot $\knot$ is both left and path
equivalent to one with ${\bf a}_d = (1, 1, \dots, 1)$.
(Thus this knot
is in ${\mathcal V}^n_d$; see \xsection
\ref{section-spaces}).

\medskip

This follows if $L$ is an orientation-preserving linear
transformation with the property that $L({\bf a}_d) = (1, 1, \dots, 1)$.
Let $(d_1, d_2, \dots, d_n)$ be the {\em vector degree} of $\knot$, where
$d_i$ is the degree of the $i$-th component of $\knot$.

\medskip

\noindent {\bf Reduction 3a}.  The knot $\knot$ is left equivalent to
a knot of vector degree $(d'_1, d'_2, \dots, d'_n)$ with $d'_1 < d'_2
< \dots < d'_n = d $ and with leading coefficients 1. 

\medskip

In fact, let $L$ be the linear transformation which reduces the
coefficient matrix to reduced row echelon form. 
This reduced form is unique.  Note, however, that the transformation
$L$ may not necessarily preserve orientation. Using an
orientation-preserving transformation gives the following result:

\medskip

\noindent {\bf Reduction 3b}.  The knot $\knot$ is both left and path
equivalent to a knot of vector degree $(d'_1, d'_2, \dots, d'_n)$ with
$d'_1 < d'_2 < \dots < d'_n = d $ and with leading coefficients
1, except for the last row, where the leading coefficient may be $\pm 1$.

\medskip

\noindent {\bf Reduction 4}.  The knot $\knot$ is left equivalent to a
knot of vector degree $(d'_1, d'_2, \dots, d'_n)$ with $d'_1 < d'_2 < \dots
d'_n \leq d $, where $d'_i$ is not in the semigroup generated by
nonnegative integral combinations of $d'_1, d'_2, \dots, d'_{i-1}$ for
$2 \leq i \leq n$.

\medskip

The transformation used here is a nonlinear shear.
This reduction is a more refined version of the one above, and it may
reduce the degree of $\knot$;
such reductions are not possible for knots in projective space.

\begin{proof}
For simplicity we prove this in the case $n = 3$; the
general case is similar. Suppose that $\knot$ has vector degree $(p,
q, r)$.  We may assume that $\knot$ is in the form of the previous
reduction.  If $ q = mp$ for a positive integer $m$, let $A:\real^3
\mapsto \real^3$ be the orientation preserving polynomial map
$$ P(x,y,z) = (x, y - x^m, z).$$
Note that $P$ has a polynomial inverse 
$Q(x', y', z') = (x',y'+{x'}^m, z')$ 
and that the leading
terms in $P \circ \knot$ have vector degree $(p, q', d)$ with $q'$ strictly
less than $q$. 
Similarly, if $d= mp + nq$ with $m,n$ nonnegative
integers, not both zero, consider the polynomial map $P :\real^3
\mapsto \real^3$ defined by
$$ P(x,y,z) = (x, y, z - x^m y^n). $$
Again, $P$ is an automorphism.  The leading terms in $P \circ \knot$ now have
degree $p, q, r$ with $r$ strictly less than $d$.
\end{proof}

These transformations are tame (see the end of \xsection
\ref{section-equivalence}. If only an even number of row switches are
involved then these knots are path equivalent. (See the remarks at the
end of \xsection \ref{section-equivalence}.)


\medskip

\subsection{Finding lower bounds on the degree}

Next we assemble some tools to examine the second question of finding
lower bounds for degree of a polynomial representation of a
given topological knot. In general these lower bounds are rather weak
and more work needs to be done here.

\subsubsection{The crossing number}
Recall that the crossing number of a topological knot $K$ (thought of
as a subset of $\real^3$) is the least number of crossings in any
planar projection of any $K'$ with the same knot type as
$K$. The ``crossing number'' of a polynomial knot refers to the crossing
number of its one-point compactification; the same applies to ``bridge
number'' and so forth.

The following two results are due to Lee Rudolph (Mount Holyoke
seminar lecture, 1995 (see \cite{Auerbach}).

\begin{lemma}
Let $\alpha(t) = (x(t), y(t))$ be a parameterized curve in the real
plane, where $x(t)$ and $y(t)$ are polynomials in $t$ of degree $d_x$
and $d_x$ respectively. Assume that $\alpha'(t) \neq {\bf 0}$ for all
$t$, and that the self-intersections of this curve are double points
(transverse intersections). Then this curve has at most
$(1/2)(d_x-1)(d_y -1)$ double points.
\end{lemma}



\begin{proof}
This curve has a double point at $(x_0,
y_0)$ if and only if there are $r \neq s$ so that $x(r) = x(s) =
x_0$ and $y(r) = y(s) = y_0$. Such $r$ and $s$ are solutions to
the equations
$$ \frac{x(r) - x(s)}{r-s} = 0$$
$$ \frac{ y(r) - y(s)}{r-s} = 0.$$
The first equation has degree $d_x - 1$, and the second
$d_y - 1$, so by Bezout's theorem these equations have at most
$(d_x -1)(d_y -1)$ intersections in the real
plane. Thus there are at most $(d_x -1)(d_y -1)$ ordered pairs
$(r,s)$ giving crossings. Since a crossing is specified by an
unordered pair, there are at most $(1/2)(d_x -1)(d_y-1)$ crossings.
\end{proof}

\begin{proposition}
\label{proposition-cross}  If $\knot(t)$ is a polynomial knot of
degree $d$ and crossing number $c$, then  $$ c \leq
(1/2)(d-2)(d-3).$$
\end{proposition}

\begin{proof}
Without loss of generality we may assume that $z(t)$ has degree
$d$. By  Reductions 3 and 4 there is a polynomial
knot $\knot_1(t) = (x_1(t), y_1(t), z(t))$ with the same knot type as
$\knot$ and $x_1(t)$ of degree at most $d-1$ and $y_1(t)$ of degree at
most $d-2$. The result then follows from the above lemma.
\end{proof}

The knot projection in the complex projective plane has exactly
$(1/2)(d-1)(d-2)$ intersections; we return to this topic in
\cite{Durfee-O'Shea-rational}. 

\subsubsection{The bridge number}
Recall that the bridge number of a topological knot $K$ can be defined
as the minimum number of local maxima of $K'$ in the direction $v$,
over all $K'$ with the same knot type as $K$, and all directions $v$.

Schubert \cite[Satz 10]{Schubert} shows that a torus
knot of type $(p,q)$ with $p < q$ has bridge number $p$, and also
(Satz 7) that the bridge number of a connected sum of knots is the sum
less one of their bridge numbers.

The following proposition is also due to Lee Rudolph (see
\cite{Auerbach}). Unfortunately it is rather weak since many knots
have bridge number two.

\begin{proposition}
\label{proposition-bridge} If $\knot$ is a polynomial knot of degree
$d$ and bridge number $b$, then
$$ b \leq (1/2)(d-1).$$
\end{proposition}

\begin{proof}
We find an upper bound for the bridge
number taking the $z$ axis as the direction.  
Let $e$ be the
degree of $z(t)$. By Reduction 3 above we may assume that $e \leq
d-2$. 
Without loss of generality we may assume that the coefficient of the
highest term in $z(t)$ is positive.
If $e$ is odd then the polynomial $z(t)$ has at most $(1/2)(e-1)$
local maxima, and if $e$ is even it has at most $(1/2)e - 1$.
Hence the spherical completion $\bar{\knot}$ of $\knot$ has at most
$(1/2)(e+1) \leq (1/2)(d-1)$ local maxima in the $z$ direction (since
it has a maximum at infinity).
\end{proof}

\subsubsection{The superbridge number}
Kuiper \cite{Kuiper} defines the superbridge number of a topological
knot $K$ as the minimum, over all $K'$ of the same knot type as $K$,
of the maximum number of local maxima of $K'$ in the direction $v$,
for all directions $v$.  He then shows that a torus knot of type
$(p,q$) has superbridge number $q$ if $p < q < 2p$, and has
superbridge number $2p$ if $2p < q$.

\begin{proposition} (REU 1998)
\label{proposition-superbridge} If $\knot$ is a polynomial knot of
degree $d$ and superbridge number $s$, then
$$s \leq (1/2)(d+1).$$
\end{proposition}

\begin{proof}
We find an upper bound for $s$ using the knot $\bar{\knot}$.  
For any rotation of $\bar{\knot}$ the degree of the $z$ coordinate is at
most $d$.
An argument as above proves
the proposition.
\end{proof}

\subsection{Applications}
\label{section-applications}
Now let examine Questions 2 and 3, proceeding degree by degree.  Let
$\knot$ be a polynomial knot.  (The unknot can be represented in all
degrees by $\knot(t) = (t^d, t, t)$.)

\bigskip

\begin{itemize}
\item $deg(\knot) \leq 4$: The crossing number (and bridge number) of
  $\knot$ is at most one, so $\knot$ is the unknot. This also follows
  since the space of such polynomial knots is path
  connected (\xsection
  \ref{section-spaces}). 
\item $deg(\knot) \leq 5$: The crossing number of $\knot$ is at most three,
hence $\knot$ is either the trefoil or the unknot.  Shastri's equations
for the trefoil are of degree five.  Equations for the mirror image of
the trefoil are obtained by changing the sign of one equation, or
switching two of them. 
\item $deg(\knot) \leq 6$: The crossing number of $\knot$ is at most 10.
McFerons' equations above for the figure-eight knot have degree
six. It seems unlikely that knots of five crossings or more can be
represented by equations of degree six, though further methods are
needed here.
\item $deg(\knot) = 7$: By Proposition \ref{proposition-superbridge}
this is the minimum degree for equations of the (2,5)-torus knot,
since it has superbridge number four.  Seven is also the minimum
degree for equations of the torus knot of type (3,4), since it has
bridge number 3. Equations for these knots are given above.
\end{itemize}

\subsection{Further methods}
Let $\knot(t) = (x(t), y(t), z(t))$ be a polynomial knot, $u(t)$ a
polynomial and $\epsilon > 0$ a suitably small number.
The REU 2006 group showed that $\tilde{\knot}(t) = (x(t), y(t), z(t) +
\epsilon u(t))$ is a polynomial knot with the same topological knot
type as $\kappa(t)$. Starting with $\kappa$ it is then possible to do
the above and then apply Reduction 4 to obtain equations of
lower degree with the same knot type as $\kappa$. 

Mui \cite{Mui-thesis} finds bounds in the cases where the polynomials
are sparse or have few monomials.

\subsection{Lexicographic ordering}
\label{lex-section}

 One could also ask Questions 2 and 3 for the vector degree of $(x(t),
y(t), z(t))$ (cf. Reduction 4 of the previous section).  Mishra
\cite{Mishra00} considers lexicographical (dictionary) order on the
degree vector of a polynomial knot, and finds the minimal vector
degree of torus knots and two-bridge knots (see the references in
\ref{section-other-examples}).

Note that a polynomial representation of minimal lexicographical order
may not be of minimal degree. For example, Shastri's representation of
the figure-eight knot, with vector degree $(3,5,7)$, has lower
lexicographical order than McFeron's representation, which has vector
degree $(4,5,6)$, but the latter has lower degree.

\subsection{Curvature}
Let $k$ be the total absolute curvature of a polynomial knot of degree
$d$.  Brutt \cite{Brutt} shows that $k \leq \pi(d-1)$; 
she uses Milnor's result \cite[Thm 3.1]{Milnor} that the average number of
local extrema of a curve in $\real^n$ is $(1/\pi)\mbox{vol}
(\sphere^{n-1})$.

\subsection{Polygonal knots}
\label{section-polygonal}
It is interesting to compare the above results for polynomial knots
with the corresponding results for polygonal (stick) knots. (A {\em
polygonal knot} is an embedded polygon in $\real^3$; this is the
classical type of knot.)
Let $e(K)$ be the number of edges in a stick
knot $K$. For the following results, see for instance \cite{Calvo,
Jin-Kim, Meissen, Randell}; the methods
of proof are in general quite different from those used for
polynomial knots. For the structure of the space of polygonal knots,
see the next section.

\begin{itemize}
\item $e(K) \leq 5$: $K$ is the unknot.
\item $e(K) = 6$: $K$ is the trefoil or the unknot.
\item $e(K) = 7$: $K$ is the figure-eight or one of the
  above.
\item $e(K) = 8$: $K$ is either a prime knot of six or fewer
  crossings, a square or granny knot, the knot $8_{19}$ (the
  (3,4)-torus knot) or the knot $8_{20}$.
\item  All knots with crossing number seven can be constructed with
  nine edges.
\end{itemize}


\section{Spaces of knots}
\label{section-spaces}


\subsection{Polynomial knots}

In this section we describe various results on the topology of 
the space of polynomial knots.
Recall that ${\mathcal K}^n_d$ is the space of polynomial knots $\real^1
\to \real^n$ of the form
$$
\begin{cases}
x_1(t) & =  a^1_d t^{d} + a^1_{d-1}t^{d-1} + \dots a^1_1 t + a^1_0\\
x_2(t) & =  a^2_d t^{d} + a^2_{d-1}t^{d-1} + \dots a^2_1 t + a^2_0\\
       & \ \ \ \ \ \ \vdots \\
x_n(t) & =  a^n_d t^{d} + a^n_{d-1}t^{d-1} + \dots a^n_1 t + a^n_0.
\end{cases}
$$
with $a^k_d \neq 0$ for some $k$.
Let ${\mathcal V}^n_d \subset {\mathcal K}^n_d$ denote the subspace of
knots of the form
$$
\begin{cases}
x_1(t) & =  t^{d} + a^1_{d-1}t^{d-1} + \dots a^1_1 t \\
x_2(t) & =  t^{d} + a^2_{d-1}t^{d-1} + \dots a^2_1 t \\
       & \ \ \ \ \ \ \vdots \\
x_n(t) & =  t^{d} + a^n_{d-1}t^{d-1} + \dots a^n_1 t
\end{cases}
$$ 
with $d \geq 1$.
Both these spaces have a finite number of path
components since they are semi-algebraic sets (see for instance
\cite[Theorem 2.2.1]{Benedetti-Risler}).
Vassiliev \cite{Vassiliev}, by analyzing the
discriminant set (the set of singular knots), proves the following
results:
\begin{itemize}
\item For $n \geq 2$ the space ${\mathcal V}^n_3$ is
contractible (see also \cite{Durst}).
\item For $n \geq 2$ the space ${\mathcal V}^n_4$ 
is homology equivalent to $\sphere^{n-2}$. 
They are homotopy equivalent if $n \geq 4$.
\item For even $d$ there is a product decomposition  ${\mathcal V}^3_d
  = X \times S^1$ for some space $X$.
\item For $n \geq 5$ and all $d$ the space ${\mathcal V}^n_d$ is
  simply-connected. 
\end{itemize}
If $d \leq 4$ these results show that $K \in {\mathcal V}^n_d$ is
unknotted. (For a simpler proof see \xsection \ref{section-applications}).

Let ${\mathcal V}^n_{d_1, d_2, \dots, d_n}$ denote the space of maps
$\real^1 \to \real^n$ of the form
$$
\begin{cases}
x_1(t) & =  t^{d_1} + a^1_{d_1-1}t^{d_1-1} + \dots a^1_1 t \\
x_2(t) & =  t^{d_2} + a^2_{d_2-1}t^{d_2-1} + \dots a^2_1 t \\
       & \ \vdots  \\
x_n(t) & =  t^{d_n} + a^n_{d_n-1}t^{d_n-1} + \dots a^n_1 t. 
\end{cases}
$$ 
Kim, Stemkoski and Yuen (REU 2000) show that this space has a conical
structure:
\begin{proposition}
Each polynomial knot $\knot$ in ${\mathcal V}^n_{d_1, d_2, \dots,
  d_n}$ lies on a curve of polynomial knots which are left-right
  equivalent to $\knot$. These curves have endpoint the map $t \mapsto
  (t^{d_1}, t^{d_2}, \dots, t^{d_n})$. 
\end{proposition}

\begin{proof}
For each $u \geq 0$, let $\knot_u$ be given by
$$
\begin{cases}
x_1(t) & =  t^{d_1} + a^1_{d_1-1}ut^{d_1-1} + \dots + a^1_1 u^{d_1-1}t \\
x_2(t) & =  t^{d_2} + a^2_{d_2-1}u t^{d_2-1} + \dots a^2_1 u^{d_2 - 1} t \\
       & \ \vdots  \\
x_n(t) & =  t^{d_n} + a^n_{d_n-1}ut^{d_n-1} + \dots + a^n_1 u^{d_n-1}t. 
\end{cases}
$$ 
For $u > 0$ we have 
$$ 
\knot_u(t) = (u^{d_1} \, x_1(t/u), \ u^{d_2} \, x_2(t/u), \dots,
u^{d_n} \, x_n(t/u)).
$$ 
Thus $\knot_1 = \knot$ and $\knot_u$ is left-right equivalent to
$\knot_1$  for $0 < u \leq 1$.
For $u = 0$ we have
$$\knot_0 (t) =  (t^{d_1}, t^{d_2}, \dots, t^{d_n}).$$

\end{proof}
Note that this proof also shows that every knot in 
${\mathcal V}^n_{d_1, d_2, \dots, d_n}$  lies
in the deformation space of the map $t \to (t^{d_1}, t^{d_2}, \dots, t^{d_n})$.

\subsection{Polynomial knots of degree five}
\label{section-gang3}

Kim, Stemkoski and Yuen \cite{gang-three} analyze the structure of the space
${\mathcal K}_5^3$ of polynomial knots of degree five in
dimension three. In this section we sketch some of their results.

We call a polynomial knot of degree five an ``S-trefoil'' if it is
topologically equivalent (in the sense of \xsection
\ref{section-equivalence}) to the Shastri trefoil of the introduction.
(In fact the Shastri trefoil is right-handed.)
A polynomial knot topologically equivalent to the mirror image of the
Shastri trefoil will be called an ``$\bar{\text S}$-trefoil''.

Let ${\mathcal{V}}_{3,4,5}$ be polynomial knots of the form

$$
\begin{cases}
x(t) &= t^3 + a_2t^2 + a_1t \\
y(t) &= t^4 + b_3t^3 + b_2t^2 + b_1t  \\
z(t) &= t^5 + c_4t^4 + +c_3t^3 + c_2t^2 + c_1t 
\end{cases}
$$
and let 
$$
{\mathcal{V}}^s_{3,4,5} \subset {\mathcal{V}}_{3,4,5} 
$$ be the knots with $a_2 = b_3 = c_3 = c_4 = 0$.  


\begin{lemma}
The region of knots in
${\mathcal{V}}^s_{3,4,5}$ whose projection is a trefoil is nonempty
and contractible.  
\end{lemma}

This lemma is proved by direct computation.  
The next lemma follows by an analysis of the knots
with these projections.

\begin{lemma}
The region of S-trefoils in ${\mathcal{V}}^s_{3,4,5}$ is nonempty and
contractible. There are no $\bar{\text S}$-trefoils in
${\mathcal{V}}^s_{3,4,5}$. 
\end{lemma}

\begin{lemma} 
\label{lemma-deformation-retraction}
There is a deformation retraction
$$ {\mathcal{V}}^s_{3,4,5} \leftarrow {\mathcal{V}}_{3,4,5}.
$$
\end{lemma}

In fact, this deformation retraction is constructed using left-right
equivalences.  
Thus the knots in ${\mathcal{V}}_{3,4,5}$ are topologically
equivalent to those in ${\mathcal{V}}^s_{3,4,5}$.

\begin{corollary}
\label{corollary-4-9}
The region of S-trefoils in ${\mathcal{V}}_{3,4,5}$ is nonempty and
contractible. There are no $\bar{\text S}$-trefoils in ${\mathcal{V}}_{3,4,5}$.
\end{corollary}

Let $\bar{{\mathcal V}}_{d_1, d_2, 5}$ be the image of ${\mathcal
    V}_{d_1, d_2, 5}$ under the map ${\bf x} \to - {\bf x}$. Since
    this map takes a trefoil to its mirror image, the following is
    also true:

\begin{corollary}
The region of $\bar{\text S}$-trefoils in $\bar{{\mathcal{V}}}_{3,4,5}$ is
nonempty and contractible. There are no S-trefoils in
$\bar{{\mathcal{V}}}_{3,4,5}$.
\end{corollary}


A principal consequence of the above results is the following
proposition and its corollary:

\begin{proposition}
The region 
of S-trefoils in ${\mathcal K_5^3}$ is nonempty and connected.
\end{proposition}

\begin{corollary} The region 
of $\bar{\text S}$-trefoils in ${\mathcal K_5^3}$ is nonempty and
connected.
\end{corollary}

The following argument proves the above proposition:
Let $\knot \in {\mathcal K_5^S}$, where ${\mathcal K_5^S} = \{ \knot
\in {\mathcal{K}_5^3} : \knot \text{ is an S-trefoil}\}$.  
By Reduction 3b (Section \ref{section-reduced-form}), $\knot$ is
both left and path equivalent to a knot $\tilde{\knot}$ in ${\mathcal
V}_{d_1, d_2, 5}$ or $\bar{{\mathcal V}}_{d_1, d_2, 5}$, where $d_1 <
d_2 < 5$.
Since left-equivalence (and path-equivalence) implies topological
equivalence (\xsection \ref{section-equivalence}), the path of knots
from $\knot$ to $\tilde{\knot}$, including the endpoints, lies in
${\mathcal K_5^S}$.  
In fact $\tilde{\knot}$ is in ${\mathcal  V}_{3,4, 5}$ or 
$\bar{{\mathcal V}}_{3,4, 5}$, since
otherwise it would be unknotted (\xsection \ref{section-applications}).
Since $\tilde{\knot}$ is an S-trefoil, Corollary \ref{corollary-4-9} implies
that it must be in ${\mathcal V}_{3,4, 5}$.
If $\knot' \in \mathcal{K}_5$ is another S-trefoil, by the
above argument it is also connected by a path of knots in ${\mathcal
  K_5^S}$ to an S-trefoil $\tilde{\knot}' \in {\mathcal{V}}_{3,4,5}$. 
Since the region ${\mathcal{V}}_{3,4,5}$ is connected (Corollary
\ref{corollary-4-9}), the proposition follows.
The corollary is a simple consequence.

\subsection{Polygonal knots}
\label{section-stick-structure}
It is interesting to compare these results with similar ones for the space
${\mathcal P}_k$ of oriented polygonal knots in $\real^3$ with $k$ edges.  
\begin{itemize}
\item ${\mathcal P}_3$, ${\mathcal P}_4$, and ${\mathcal P}_5$ are
  connected. (This result is attributed to Kuiper.)
\item ${\mathcal P}_6$ has 5 path components: one of unknots, and two each of
  right-handed and left-handed trefoils \cite{Calvo}.
\item ${\mathcal P}_7$ has 5 path components: unknots, right-handed trefoils,
  left-handed trefoils, and two of figure-eight knots \cite{Calvo}.
\item ${\mathcal P}_8$ has at least 20 path components \cite{Calvo}.
\end{itemize}


\section{Complex polynomial knots}
\label{section-complex}

Let $\alpha: \complex^1 \to \complex^n$ be a polynomial map; in
particular $\alpha$ can be the complexification of a real polynomial
map.  (Maps of this type in the projective case have been studied by
Viro \cite{Viro}.)  We say that $\alpha$ is an embedding if it is an
injection and $\alpha'(t) \neq 0$ for all $t \in \complex$; in case we
call $\alpha$ a {\em complex polynomial knot}.  
Note that the set of polynomial knots of degree $d$ such that its
complexification is an embedding is dense in ${\mathcal K}^n_d$.

For $n \geq 3$ complex polynomial knots are topologically unknotted.
Furthermore the parameter space of such knots of degree $d$ is connected
since the subset of maps with singularities has codimension two. 

The question of left-right equivalence is complicated.  A
polynomial embedding $\alpha: \complex \to \complex^n$ is {\em
rectifiable} if there is a polynomial automorphism $P$ of $\complex^n$
such that $P \circ \alpha = \iota$, where $\iota(t) = (t, 0, \cdots,
0)$. In other words, it is rectifiable if it is left equivalent to a
linear map.  A basic question is whether every polynomial embedding is
rectifiable.  This is true for $n = 2$ by the theorem of
Abhyankar-Moh, a surprising result since topological knotting is
possible in this dimension.  It is also true for $n \geq 4$
\cite{Craighero, Jelonek}.  For $n = 3$ this is apparently unknown.
Shastri's work initiated in this question; he conjectured that the
trefoil knot of the introduction is not rectifiable.  For this circle
of ideas see for example \cite{van-den-Essen}.

Next we describe some connections with commutative algebra and the
standard correspondences of algebraic geometry.  If $\alpha =
(\alpha_1, \dots, \alpha_n): \complex \to \complex^n$ is a polynomial
map, let ${\hat \alpha}: \complex[x_1, \dots, x_n] \to \complex[t]$ be
defined by ${\hat \alpha}(x_i) = \alpha_i(t)$.
%
%
%
%
Conversely, if $\beta: \complex[x_1, \dots, x_n] \to \complex[t]$ is a
map of rings,
then ${\hat \beta} = ({\hat \beta_1}, \dots, {\hat \beta_n})$
is defined by ${\hat \beta _i}(t) = \beta(x_i)$.  
This provides a 1-1 correspondence between polynomial maps
$\complex \to \complex^n$ and 
maps of rings $\complex[x_1, x_2, \dots, x_n] \to \complex[t]$.


\begin{proposition}
The map $\alpha$ is an embedding if and only if the map ${\hat \alpha}$ is
surjective.
\end{proposition}

This proposition makes Gr\"{o}bner basis methods useful for
determining whether a complex map is an embedding. 
(See also \cite{Siegel}.)

\begin{proof}
Let $V \subset \complex^n$ denote the image of $\alpha$, and let $I
\subset \complex[x_1, \dots x_n]$ be its ideal.
Let 
$\alpha^*: \complex[x_1, \dots x_n]/I \to \complex[t]$ be the map dual
to $\alpha$.
Proposition A.2.12 (p. 417) of \cite{Greuel-Pfister} asserts that
$\alpha^*$ is surjective exactly when $V \subset \complex^n$ is a
closed subvariety and $\alpha: \complex \to V$ is an isomorphism.

 Since the map $\hat{\alpha}$ factors through $\alpha^*$, the former
is surjective exactly when the latter is.  Also $\alpha$ is an
isomorphism exactly when it is bijective and its derivative $D_t \,
\alpha$ is nonsingular for all $t \in \complex$.  
Finally, $V$ is always a closed subvariety since $\alpha$ is a
polynomial map.
\end{proof}

Here is a simpler proof of the forwards implication:
Let $\alpha$ be as above,
and $t \in \complex^1$.
Since ${\hat \alpha}$ is surjective,
there is a polynomial $p \in \complex[X_1, \dots, X_n]$ such that 
$p(\alpha(t)) = t$. Taking derivative with respect to $t$ gives
$(Dp) \circ \alpha' = 1$.  Thus $\alpha'(t) \neq 0$ for all $t$.  Next
we show that $\alpha$ is injective.  Let $s,t \in \complex^1$ with the
property that $\alpha(s) = \alpha(t)$.  Applying $p$ to both sides
gives $s = t$. 

\begin{example}
(1) The Shastri trefoil extended to $\complex$ is
an embedding since $yz - x^3 -5xy + 2z - 7x$ maps to $t$. Similarly
for the complexified figure-eight knot the polynomial
$x^2z-xy^2-7x^2y-23x^3-3z+22y+71x$ maps to $t$. (See \cite{Shastri}.)

\noindent (2). The cusp $\alpha(t) = (t^2, t^3)$ is not an embedding;
the map is defined by $x \mapsto t^2$ and $y \mapsto t^3$ is not
surjective. The same is true of the double point $\alpha(t) = (t^2,
t^3-1)$.
Also the proposition is not true over $\real$; for example
$\alpha(t) = t^3 + t$ is an embedding but $\hat{\alpha}$ is not
surjective. (Jason Starr.)
\end{example}


\bigskip

Papers by Mount Holyoke REU students cited in the references below can
be found at www.mtholyoke.edu/acad/math/reu.

\end{document}